\documentclass[aps,pre,twocolumn,superscriptaddress,nofootinbib]{revtex4-2}

\usepackage{amsmath,amssymb,amsfonts,bm}
\usepackage{graphicx}
\usepackage{textcase}
\usepackage{xcolor}
\usepackage{hyperref}
\usepackage{tikz}
\usepackage{pgfplots}

\begin{document}

\title{On the Second Moment of the Maximum of Fractional Brownian Motion}

\author{Safari Mukeru and Mmboniseni P. Mulaudzi}
\email{mukers@unisa.ac.za; mulaump@unisa.ac.za}
\affiliation{Department of Decision Sciences, School of Economic Sciences, University of South Africa, P.O. Box 392, Pretoria 0003, South Africa}

\date{\today}

\begin{abstract} Let \(B_H\) be standard fractional Brownian motion and \(M_H=\sup_{0\leq t\leq1}B_H(t)\). Motivated by Spitzer's identities for the maximum of partial sums with exchangeable increments, we propose a one-dimensional integral approximation for \(\langle M_H^2\rangle\). The expression is exact at \(H=1/2\) and \(H=1\), agrees closely with numerical simulations for intermediate \(H\), and yields a correction near \(H=1/2\) close to that obtained by perturbative methods. We also briefly discuss the limiting behavior as \(H\downarrow0\). The unexpected accuracy of an exchangeability-based construction for a process with generally nonexchangeable increments points to a structural connection that remains to be understood. It also suggests that identities derived under exchangeability may provide a useful source of candidate approximations for extremal quantities of fractional Brownian motion when exact methods are unavailable. \end{abstract}

\maketitle

\section{Introduction}

The maximum of a stochastic process is a fundamental quantity in problems involving extreme excursions, threshold exceedances, and risk over a finite observation period. When the underlying dynamics exhibit long-range dependence, however, the distribution and moments of the maximum are generally difficult to determine.
Fractional Brownian motion (fBm) provides a canonical Gaussian framework for studying this problem because it combines self-similarity with stationary, temporally correlated increments. 

Applications involving correlated fluctuations arise in several areas. In hydrology, fractional Brownian motion and fractional Gaussian noise have been used to describe geophysical time series and spatial heterogeneity in subsurface systems \cite{MandelbrotVanNess1968, MolzLiuSzulga1997}. 
In such settings, information about the maximum is relevant to the assessment of extreme water levels, reservoir inflows, and other unusually large hydrological excursions. 
Related models occur in anomalous transport, where temporally correlated motion may describe the effective trajectory of a tagged particle in a crowded or viscoelastic environment \cite{MetzlerJeonCherstvyBarkai2014,Sokolov2012}. In biological applications, fBm-based models have been considered for the motion of passive tracers and molecules in the cytoplasm and other crowded intracellular environments \cite{RunfolaVitaliPagnini2022}.
Fractional Brownian motion has also been used in finance to represent persistent or antipersistent fluctuations, particularly in models of volatility with long memory \cite{ComteRenault1998}. These applications motivate the study not only of the expected maximum, but also of its variability. 
We consider standard fractional Brownian motion \[ \{B_H(t):t\geq0\}, \] defined as the centered Gaussian process with Hurst exponent $H\in(0,1)$, $B_H(0)=0$, and covariance \begin{equation} \left\langle B_H(t)B_H(s)\right\rangle = \frac12 \left( t^{2H}+s^{2H}-|t-s|^{2H} \right). \label{eq:covariance} \end{equation} 
We also consider its continuous extension to $H=1$, corresponding to the degenerate Gaussian process $B_1(t)=tZ$, where $Z\sim N(0,1)$. For $H=1/2$, the process reduces to standard Brownian motion. For $H\neq1/2$, it is non-Markovian, with negatively correlated increments when $H<1/2$ and positively correlated increments when $H>1/2$ \cite{MandelbrotVanNess1968,Nourdin2012}. Let \begin{equation} M_H(T) = \sup_{0\leq t\leq T}B_H(t) \end{equation} denote the running maximum over $[0,T]$. 
Self-similarity gives \begin{equation} M_H(T)\overset{d}{=}T^H M_H \mbox{  with  }
M_H := M_H(1)  \end{equation}
and consequently \begin{equation} \left\langle M_H(T)^q\right\rangle = T^{qH}\left\langle M_H^q\right\rangle, \qquad q>0, \label{eq:moment-scaling} \end{equation} whenever the corresponding moment exists. The analysis may therefore be reduced to the unit interval.

In \cite{MukeruMulaudziFirstMoment}, we studied the first moment of the maximum and proposed the approximation \begin{equation} \left\langle M_H\right\rangle \approx \left\langle M_H\right\rangle_{\mathrm{app}} := \frac{1}{H\sqrt{2\pi}}. \label{eq:first-moment-approximation} \end{equation}

The present paper continues that line of investigation by considering the second moment \[ \left\langle M_H^2\right\rangle \] and the variance, \begin{equation} \operatorname{Var}\!\left(M_H\right) = \left\langle M_H^2\right\rangle - \left\langle M_H\right\rangle^2. \label{eq:variance-definition} \end{equation} Accurate expressions for both quantities provide a more informative statistical description of the maximum than its expectation alone and may facilitate moment-based comparisons with empirical or simulated data. Despite the Gaussian structure of fBm, its non-Markovian dependence makes its extreme-value statistics considerably more difficult to determine than those of ordinary Brownian motion. Previous investigations have established persistence and small-value asymptotics and developed perturbative expansions for the distribution of the maximum and related path functionals \cite{Molchan1999,WieseMajumdarRosso2011,DelormeWiese2015,DelormeWiese2016,SadhuDelormeWiese2018}. These results provide detailed local information, particularly near \(H=1/2\), but no elementary closed-form expression for the second moment is presently known for general \(H\). The moments of \(M_H\) also admit a direct physical interpretation in anomalous transport. If the displacement of a tracer is modeled by \[ X(t)=\sqrt{2D_H}\,B_H(t), \] where \(D_H\) is the generalized diffusion coefficient, then \[ \langle X(t)^2\rangle=2D_Ht^{2H}. \] Self-similarity implies that the largest displacement attained over an observation time \(T\) satisfies \[ X_{\max}(T) \overset{\mathrm d}{=} \sqrt{2D_H}\,T^H M_H, \] and hence \[ \operatorname{Var}\!\left[X_{\max}(T)\right] = 2D_H T^{2H}\operatorname{Var}(M_H). \] Thus, the moments of \(M_H\) characterize the magnitude and variability of the furthest excursion, which are relevant to boundary- and target-reaching problems. This interpretation applies, for example, to tagged-particle motion in single-file systems when the effective dynamics are described by fBm \cite{LimTeo2009,AkintundeEtAl2025}. Although the present work continues our earlier investigation of the mean maximum, the second-moment problem is not a direct extension of the first-moment calculation. Motivated by the second-moment identity arising from Spitzer's combinatorial analysis of partial sums with exchangeable increments, that is, increments whose joint distribution is unchanged under permutation, we propose a one-dimensional integral approximation for \(\langle M_H^2\rangle\). The second-moment calculation involves a distinct exchangeability identity and, in the continuum formulation, a nontrivial contribution from the joint positive parts of two correlated fBm increments. Its extension to fBm is necessarily heuristic, since fractional Gaussian noise is not exchangeable except at the Brownian point and at the degenerate endpoint \(H=1\). Nevertheless, the proposed expression reproduces the exact second moment at these two solvable values and agrees closely with numerical simulations over the range of \(H\) investigated. We further expand the approximation about \(H=1/2\) and compare the resulting coefficient with the second-moment correction formally extracted from the first-order maximum-density expansion of Delorme and Wiese \cite{DelormeWiese2015,DelormeWiese2016}. The two coefficients are different, although numerically close. This discrepancy shows that the proposed expression should be regarded as a highly accurate approximation rather than an exact identity. The paper is organized as follows. We first recall the second-moment identity for maxima of partial sums of exchangeable random variables and formulate its continuum analogue, from which the proposed approximation is obtained. We then verify the approximation at the exactly solvable values \(H=1/2\) and \(H=1\). Next, we expand it about the Brownian point and compare the resulting coefficient with that formally extracted from the Delorme--Wiese perturbative density. Finally, we compare the approximation with numerical estimates and discuss the limiting regime \(H\downarrow0\).

\section{The second-moment approximation}
Our main approximation for the second moment is
\begin{equation} \label{eq:secon_order_appr} \left\langle M_H^2 \right\rangle \approx \left\langle M_H^2 \right\rangle_{\mathrm{app}} := \frac{1}{4H} + \frac{I(H)}{4\pi H}. \end{equation}
where
\begin{equation}
 \begin{aligned} 
 I(H) = 
 \int_0^1
 \frac{\sqrt{1-\rho_H(x)^2}
 +\rho_H(x)
 \left(\frac\pi2+\arcsin\rho_H(x)\right)}{x^{1-H}(1-x)^{1-H}}\,
 dx,
 \end{aligned}
 \label{eq:conjecture}
\end{equation}
with
\[
 \rho_H(x)=
 \frac{1-x^{2H}-(1-x)^{2H}}
 {2x^H(1-x)^H}.
\]
The quantity \(\rho_H(x)\) has a direct interpretation as the correlation between the two adjacent increments \(B_H(x)\) and \(B_H(1)-B_H(x)\): \begin{equation} \rho_H(x) = \frac{ \left\langle B_H(x) (B_H(1)-B_H(x))\right\rangle }{ \sqrt{\left\langle B_H(x)^2\right\rangle} \sqrt{\langle(B_H(1)-B_H(x))^2\rangle} }. \label{eq:rho-interpretation} \end{equation}
Thus the approximation depends only on the correlation between two adjacent fBm blocks of relative lengths $x$ and $1-x$.

Combining the separate approximations for the first and second moments gives the formal variance estimate \[ V_{\mathrm{app}}(H) = \frac{1}{4H} +\frac{I(H)}{4\pi H} -\frac{1}{2\pi H^2}. \] Although \(V_{\mathrm{app}}(H)\) remains nonnegative throughout the numerically investigated range \(H\geq1/4\), it becomes negative for sufficiently small \(H\). Indeed, since \[ I(H)\sim\frac{3\sqrt{3}}{4H}, \qquad H\downarrow0, \] we obtain \[ V_{\mathrm{app}}(H) \sim \frac{3\sqrt{3}-8}{16\pi H^2}<0. \] This breakdown is consistent with upper bounds obtained by Borovkov et al.\ \cite{BorovkovEtAl2017,BorovkovEtAl2018}, which show that \(\langle M_H\rangle\) grows at most on the order of \(H^{-1/2}\) as \(H\downarrow0\), whereas the first-moment approximation \(1/(H\sqrt{2\pi})\) grows as \(H^{-1}\). Since the supporting simulations were restricted to \(H\geq1/4\), the first-moment approximation and the derived variance estimate should be used only within this range.‌

 The limiting case \(H\downarrow0\) will be discussed separately in Section~\ref{subsec:discrete-H-zero}.



The purpose of this paper is to propose and assess the approximation in Eq.~\eqref{eq:conjecture}. The approximation is exact for Brownian motion, corresponding to $H=\tfrac12$, and at the degenerate Gaussian endpoint $H=1$. It is also strongly supported by numerical simulations over a broad range of Hurst exponents. However, its first-order expansion about $H=\tfrac12$ differs slightly from the exact first-order result obtained using the perturbative maximum-density formula of Delorme and Wiese \cite{DelormeWiese2015,DelormeWiese2016}.

\section{Partial sums of exchangeable random variables and the second moment of their maximum}
A finite collection of random variables is exchangeable if their ordering carries no statistical information: permuting the variables does not change their joint distribution. Independent and identically distributed random variables are exchangeable, although exchangeable random variables need not be independent.
Let $Y_1,\ldots,Y_n$ be $n$ exchangeable random variables, define
\[
 S_r=Y_1+\cdots+Y_r,
 \qquad
 S_r^+=\max(0,S_r),
\]
and put
\[
 \mathcal M_n=\max(0,S_1,\ldots,S_n).
\]
It is proven in \cite{AnisGharib1980} that 
\begin{equation}
 \begin{aligned}
 \langle \mathcal M_n^2 \rangle
 &=\sum_{k=1}^{n}\frac{1}{k} \langle(S_k^+)^2 \rangle
 \\
 &\quad+\sum_{r=2}^{n}
 \sum_{k=1}^{r-1} \frac{\langle S_k^+(S_r-S_k)^+\rangle}{k(r-k)}
 \end{aligned}
 \label{eq:AG-second-moment}
\end{equation}
For independent increments, $S_k$ and $S_r-S_k$ are independent and the cross moment factorizes. For exchangeable but dependent increments, the joint positive-part expectation must be retained.

Formula \eqref{eq:AG-second-moment} is based on the exchangeability property of the involved random variables. Fractional Gaussian noise is stationary but not exchangeable for general $H$, because its covariance depends on the lag $|i-j|$. However, we shall substitute the sums of the fBm blocks into \eqref{eq:AG-second-moment}, pass formally to the continuum, and study the resulting expression as a conjectural formula. No claim of finite-$n$ validity for fGn is made.

\section{Gaussian positive parts}
We require a standard bivariate Gaussian identity.
\label{lem:positive-parts}
Let $(U,V)$ be jointly Gaussian random variables with 
\[\langle U \rangle = \langle V \rangle  = 0, \quad 
\sqrt{\langle U^2 \rangle} = \sigma_U, \sqrt{\langle V^2 \rangle} = \sigma_V\] and  
 $$\frac{\langle U V \rangle}{\sigma_U \sigma_V} = \rho \in [-1,1]$$
 Then
\begin{equation}
 \langle U^+V^+\rangle 
 =\frac{\sigma_U\sigma_V}{2\pi}
 \left[
 \sqrt{1-\rho^2}
 +\rho\left(\frac\pi2+\arcsin\rho\right)
 \right].
 \label{eq:positive-parts}
\end{equation}
For convenience, set 
\begin{equation}
 \mathcal G(\rho)
 :=\sqrt{1-\rho^2}
 +\rho\left(\frac\pi2+\arcsin\rho\right).
 \label{eq:Gplus}
\end{equation} and then 
 \begin{equation} 
 \langle U^+V^+\rangle 
 =\frac{\sigma_U\sigma_V}{2\pi} \mathcal{G}(\rho)
\end{equation}
\section{Formal continuation from exchangeable sums to fBm}
For a fixed $n$, set
\begin{eqnarray*}
 X_i=B_H(i/n)-B_H((i-1)/n), \\
 S_r=\sum_{i=1}^rX_i=B_H(r/n).
\end{eqnarray*}
We now insert these variables formally into \eqref{eq:AG-second-moment}.

\subsection{The one-block term}
Since $S_r\sim N(0,(r/n)^{2H})$,
\[
 \langle (S_r^+)^2 \rangle =\frac12\left(\frac rn\right)^{2H}.
\]
Hence
\begin{equation}
 A_{n,H}:=\sum_{r=1}^n\frac1r\langle (S_r^+)^2 \rangle 
 =\frac1{2n^{2H}}\sum_{r=1}^nr^{2H-1}.
 \label{eq:AnH}
\end{equation}
The power-sum asymptotic gives
\begin{eqnarray}
 \lim_{n\to\infty}A_{n,H} &= &\lim_{n\to \infty}\frac{1}{2n} \sum_{r=1}^n (r/n)^{2H-1} \nonumber \\
 & = & \frac{1}{2}\int_0^1 x^{2H-1} dx =  
 \frac1{4H}.
 \label{eq:one-block-limit}
\end{eqnarray}

\subsection{The two-block term}
Fix $1\leq k<r\leq n$ and put
\[
 U=S_k=B_H(k/n)\] and 
\[ V=S_r-S_k=B_H(r/n)-B_H(k/n).\]
Their standard deviations are
\[
 \sigma_U=(k/n)^H,
 \qquad
 \sigma_V=((r-k)/n)^H.
\]
Their covariance is
\[
 \langle U V \rangle 
 =\frac1{2n^{2H}}
 \left(r^{2H}-k^{2H}-(r-k)^{2H}\right).
\]
Therefore their correlation is
\begin{equation}
 \rho_H(k,r)
 =\frac{r^{2H}-k^{2H}-(r-k)^{2H}}
 {2k^H(r-k)^H}.
 \label{eq:finite-block-correlation1}
\end{equation}
Formula \eqref{eq:positive-parts} yields
\begin{equation}
 \langle S_k^+(S_r-S_k)^+ \rangle 
 =\frac{k^H(r-k)^H}{2\pi n^{2H}} 
 \mathcal G\bigl(\rho_H(k,r)\bigr).
 \label{eq:block-positive-parts}
\end{equation}
The second term in \eqref{eq:AG-second-moment} becomes  
\begin{equation}
 D_{n,H}
 =\frac1{2 \pi n^{2H}}
 \sum_{r=2}^{n}
 \sum_{k=1}^{r-1} \frac{k^H(r-k)^H}{k(r-k)} \mathcal G\bigl(\rho_H(k,r)\bigr).
 \label{eq:BnH}
\end{equation}

Set $x=k/r$. Dividing both the numerator and denominator in \eqref{eq:finite-block-correlation1} by $r^H$ yields
\[\rho_H(k,r)=\rho_H(k/r)\] with
\begin{equation}
 \rho_H(x)
 :=\frac{1-x^{2H}-(1-x)^{2H}}
 {2x^H(1-x)^H}.
 \label{eq:rhoH}
\end{equation}
Moreover, \[ \frac{k^H(r-k)^H}{k(r-k)} = r^{2H-2} \left(\frac{k}{r}\right)^{H-1} \left(1-\frac{k}{r}\right)^{H-1}. \] Hence \[ D_{n,H} = \frac{1}{2\pi n^{2H}} \sum_{r=2}^{n} r^{2H-1} \left[ \frac1r\sum_{k=1}^{r-1} f_H\left(\frac{k}{r}\right) \right], \] where \[ f_H(x) = x^{H-1}(1-x)^{H-1}\mathcal G(\rho_H(x)). \]
Thus the inner sum in \eqref{eq:BnH} is  a Riemann sum and hence it converges to 
  \[I(H) = \int_0^1 f_H(x) dx.\]
Thus for $n \to \infty$, 
\[D_{n,H} \sim \left(\frac{1}{2\pi n^{2H}} \sum_{r=2}^n r^{2H - 1} \right) I(H).\]
Consequently,
\begin{equation}
 \lim_{n\to\infty}D_{n,H}=\frac{I(H)}{4\pi H}.
 \label{eq:two-block-limit}
\end{equation}
This is how approximation Eq.~\eqref{eq:secon_order_appr} was obtained. It can also be written as 
\begin{equation} \begin{aligned}\left\langle M_H^2 \right\rangle_{\mathrm{app}} = {}& \int_0^1 \frac{\langle \left(B_H(t)^+\right)^2\rangle }{t}\,dt \\ &+ \frac{1}{2H} \int_0^1 \frac{\langle B_H(t)^+ \left(B_H(1)-B_H(t)\right)^+ \rangle }{ t(1-t) }\,dt . \end{aligned} \label{eq:second-moment-two-block} \end{equation}

\section{Exact checks}

\subsection{Brownian motion, $H=1/2$}

For $H=1/2$, the process $B_H$ coincides with standard Brownian motion. According to the reflection principle, reflecting a Brownian path about the line $y=a$ after its first hitting time of the level $a>0$ preserves its probability law. Consequently, the running maximum over the unit interval has the same distribution as the absolute value of a standard normal random variable. Hence, \begin{equation}
M_{1/2}\overset{d}{=}|Z|, \qquad  Z\sim N(0,1).
\label{eq:abdvalue_Gauss}\end{equation} 
A standard treatment of the reflection principle and the distribution of the running maximum may be found in Ref.~\cite{KaratzasShreve1991}.
When $H=1/2$, the adjacent blocks $S_k$ and $S_r - S_k$ (for $k < r$) are independent and $\rho_{1/2}(x)=0$. Hence
\[
 I(1/2)=\int_0^1x^{-1/2}(1-x)^{-1/2}dx
 =\pi.
\]
Eq.~\eqref{eq:secon_order_appr} gives
\[
 \langle M_{1/2}^2 \rangle =  \langle M_{1/2}^2 \rangle_{\mathrm{app}} = 1
\]
in agreement with Eq.~\eqref{eq:abdvalue_Gauss}.

\subsection{Rank-one endpoint}
When $H=1$, the adjacent blocks are perfectly correlated: $\rho_1(x)=1$ and $\arcsin 1 = \pi/2$. Hence
\[
 I(1)=\pi,
\]
and
\[
 \langle M_1^2\rangle =\frac14+\frac{1}{4\pi} \pi=\frac12,
\]
in agreement with $M_1=Z^+$.

\section{First-order expansion around the Brownian point}
Set $H=1/2+\varepsilon$, with $|\varepsilon|$ small. Differentiation of \eqref{eq:secon_order_appr} under the integral sign yields the first-order expansion \begin{equation} \label{eq:conjecture-first-order} \left\langle M_{1/2+\varepsilon}^2\right\rangle_{\mathrm{app}} = 1+ \left( \frac{\pi^2}{12}-2\log 2-2 \right)\varepsilon + O(\varepsilon^2). \end{equation}

We now compare the coefficient in \eqref{eq:conjecture-first-order} with the first-order correction obtained from the perturbative density of $M_H$
\begin{equation}
f_H(y)=f_0(y)\left[1+\varepsilon\bigl(G(y)+c\bigr)\right]
 +O(\varepsilon^2),
 \label{eq:DW-density}
 \end{equation}
 with
 \[ f_0(y)=\sqrt{\frac2\pi}e^{-y^2/2}.\]
The explicit correction is
\begin{equation}
 \begin{aligned}
 G(y)
 &=\frac{y^4}{6}\,{}_2F_2\left(1,1;\frac52,3;\frac{y^2}{2}\right)
 -3y^2
 \\
 &\quad+\pi(1-y^2)\operatorname{erfi}\left(\frac y{\sqrt2}\right)
 +\sqrt{2\pi}\,e^{y^2/2}y
 \\
 &\quad+(y^2-2)\left[\log(2y^2)+\gamma_E\right].
 \end{aligned}
 \label{eq:DW-G}
\end{equation}
(See \cite{WieseMajumdarRosso2011,DelormeWiese2015, DelormeWiese2016}.)
The constant $c$ ensures normalization and its value is 
\begin{equation}
 c=-\int_0^\infty f_0(y)G(y)dy = -2;
 \label{eq:c-normalization}
\end{equation}
as previously obtained in \cite{MukeruMulaudziFirstMoment}.
From \eqref{eq:DW-density},
\begin{align}
 \langle M_{1/2+\varepsilon}^2 \rangle 
 &=\int_0^\infty y^2f_H(y)dy
 \nonumber\\
 &=1+\varepsilon\int_0^\infty y^2f_0(y)(G(y)+c)dy
 +O(\varepsilon^2).
 \label{eq:DW-second-moment-start}
\end{align}
Using $c=-2$ and $\int y^2f_0=1$, we obtain
\begin{equation}
 \left.\frac d{dH}\langle M_H^2 \rangle \right|_{H=1/2}
 =J-2,
 \end{equation} with 
 \begin{equation} 
 J:=\int_0^\infty y^2f_0(y)G(y)dy.
 \label{eq:J2-definition}
\end{equation}
The integral $J$ is calculated in Appendix A by regularisation and its value is
 \[J = -4 \log2 +2.\] 
 It follows that 
 \begin{equation}
 \left.\frac d{dH}\langle M_H^2 \rangle \right|_{H=1/2}
 = - 4 \log 2.
 \label{eq:DW-second-moment-coefficient}
\end{equation}
Consequently,
\begin{equation}
 \langle M_{1/2+\varepsilon}^2 \rangle
 =1 - 4 \varepsilon \log2  +O(\varepsilon^2).
 \label{eq:DW-second-moment-expansion}
\end{equation}
The approximation predicts the first-order coefficient \[ \frac{\pi^2}{12}-2\log 2-2 \approx -2.56383, \] compared with the exact value \[ -4\log 2\approx -2.77259. \] The relative discrepancy is approximately $7.5\%$. Thus, although the approximation is not exact, it correctly captures both the sign and the principal magnitude of the first-order dependence near the Brownian point.

\section{Numerical evidence}

 We simulated fBm on an equally spaced grid and approximated \[ M_H=\sup_{0\leq t\leq1}B_H(t) \] by \[ M_{H,n}=\max_{0\leq k\leq n}B_H(k/n). \] Continuity implies that $M_{H,n}\leq M_H$ pathwise and $M_{H,n}\to M_H$ as $n\to\infty$. Thus, discretization introduces a downward bias, although finite-sample variability may cause the Monte Carlo estimate to lie above or below the exact moment. Finer grids are required for smaller $H$, where the paths are rougher.

For $H<0.25$, the numerical estimates become increasingly unstable despite the use of very fine grids. We therefore restrict our numerical analysis to $H\geq0.25$. 

We used the Davies--Harte algorithm and the Hosking algorithm, which is based on the Durbin--Levinson recursion,  to generate the fBm paths. Detailed descriptions of both algorithms are given by Dieker \cite[Secs.~2.1.1 and 2.1.3]{Dieker2004}.
The Davies--Harte algorithm is considerably faster. For the present implementation and grid sizes, the Davies--Harte circulant embedding became poorly conditioned near $H=1$. We therefore used the Davies--Harte algorithm for $0.25\leq H\leq0.92$ and Hosking's algorithm for $0.93\leq H\leq0.99$. For each value of $H$, we generated $1000$ paths. With the Davies--Harte algorithm, the grid spacing was $2^{-23}$ for $H=0.25,0.26,\ldots,0.29$, $2^{-20}$ for $H=0.30,0.31,\ldots,0.39$, and $2^{-18}$ for $H=0.40,0.41,\ldots,0.92$. For $H=0.93,0.94,\ldots,0.99$, we used Hosking's algorithm with grid spacing $2^{-15}$.
   
\begin{figure}[t]
\centering 
\begin{tikzpicture}
\begin{axis}[ width=0.50\textwidth, height=0.50\textwidth, xlabel={Hurst exponent $H$}, ylabel={Second moment }, xmin=0.25, xmax=1, ymin=0.35, ymax=2.90, major grid style={gray!30}, minor grid style={gray!15}, minor x tick num=10, minor y tick num=10,
legend style={ draw=none, fill=none, font=\small, at={(0.97,0.97)}, anchor=north east }, grid=major, grid style={gray!20}, tick label style={font=\small}, label style={font=\small} ]

\addplot[ blue, very thick, mark=none, smooth, tension=0.5] coordinates{
(0.25,2.71863) (0.26,2.55623) (0.27,2.41028) (0.28,2.27858) (0.29,2.15928) (0.30,2.05082) (0.31,1.95189) (0.32,1.86138) (0.33,1.77831) (0.34,1.70188) (0.35,1.63137) (0.36,1.56617) (0.37,1.50575) (0.38,1.44962) (0.39,1.39739)
(0.40,1.34869) (0.41,1.30319) (0.42,1.26062) (0.43,1.22073) (0.44,1.18327) (0.45,1.14807) (0.46,1.11493) (0.47,1.08369) (0.48,1.05421) (0.49,1.02635) (0.50,1.00000) (0.51,0.97504) (0.52,0.95138) (0.53,0.92893) (0.54,0.90760)
(0.55,0.88733) (0.56,0.86803) (0.57,0.84966) (0.58,0.83215) (0.59,0.81544) (0.60,0.79950) (0.61,0.78428) (0.62,0.76972) (0.63,0.75581) (0.64,0.74249) (0.65,0.72974) (0.66,0.71753) (0.67,0.70583) (0.68,0.69461) (0.69,0.68384)
(0.70,0.67352) (0.71,0.66360) (0.72,0.65408) (0.73,0.64494) (0.74,0.63616) (0.75,0.62772) (0.76,0.61960) (0.77,0.61181) (0.78,0.60431) (0.79,0.59710) (0.80,0.59017) (0.81,0.58351) (0.82,0.57711) (0.83,0.57095) (0.84,0.56504)
(0.85,0.55936) (0.86,0.55391) (0.87,0.54868) (0.88,0.54366) (0.89,0.53885) (0.90,0.53426) (0.91,0.52986) (0.92,0.52567) (0.93,0.52168) (0.94,0.51789) (0.95,0.51431) (0.96,0.51093) (0.97,0.50778) (0.98,0.50487) (0.99,0.50224)
};
\addlegendentry{Approximation values: $\langle M_H^2\rangle_{\mathrm{app}}$. }
 
\addplot[ only marks, mark=*, mark size=1.2pt, red ] coordinates {(0.25,2.737469950238958)  
(0.26, 2.589743625647868)  
(0.27, 2.453892595655813) (0.28,2.328818789182720)  (0.29,2.213522212056118)

(0.30,2.079005967465044) (0.31,1.985696072886134) (0.32,1.899028957264725) 
(0.33,1.818392273761976)  (0.34,1.743193247483318)  (0.35,1.673036104204241) 
(0.36,1.607503722489931)  (0.37,1.546229197222635) (0.38,1.488901959773556) 
(0.39,1.435188365092895)

(0.40,1.378169831996492) 

(0.41,1.330830777069768)

(0.42,1.286303579463445) 

(0.43,1.244413665958707) 

(0.44,1.204926838610628) (0.45,1.167684328218164) (0.46,1.132528167458470) (0.47,1.099391386244649) 

(0.48,1.068122284301974) (0.49,1.038541811457360)
(0.50,1.010557236944815) (0.51,0.984056810742762) (0.52,0.958944922502778)

(0.53,0.9352088479122684) (0.54,0.912675286196706) (0.55,0.891255215695312) 
(0.56,0.870866859940701) (0.57,0.851423967221645) (0.58,0.832930263412624) (0.59,0.815317688164550) (0.60,0.798478757643938) (0.61,0.782422818492819) (0.62,0.767119335010055) (0.63,0.752506143599208) (0.64,0.738525962540108) (0.65,0.725172492093191) (0.66,0.712472954454239) (0.67,0.700362024112030) (0.68,0.688783892973722) (0.69,0.677737344556548) (0.70,0.667233159280423) (0.71,0.657214767236231) (0.72,0.647667224309071) (0.73,0.6385599014330) (0.74,0.629862308449711) (0.75,0.621525372655100) (0.76,0.613535265385587) (0.77,0.605935217761880) (0.78,0.598709144192387) (0.79,0.591801858985135) (0.80,0.585222987167519) (0.81,0.578963103053756) (0.82,0.573027485446096) (0.83,0.567358032837263) (0.84,0.561945540457244) (0.85,0.556786734414239) (0.86,0.551839668712200) (0.87,0.547117547364133) (0.88,0.542617130281259) (0.89,0.538351018294072) (0.90,0.534347791856686) (0.91,0.530577723509792) (0.92,0.527066112513134) (0.93,0.515553979847224) 
(0.94, 0.512805398005332) 
(0.95,0.510612373872962) (0.96,0.509027563448415) (0.97,0.508132552644935) (0.98,0.508502672699994) (0.99,0.510943871360614) (1, 0.5) }; \addlegendentry{Simulated second moments} \end{axis} \end{tikzpicture} \caption{ Simulated second moments of the maximum of fractional Brownian motion on $[0,1]$ compared with the approximation values $\langle M_H^2\rangle_{\mathrm{app}}$. } \label{fig:simulation-theory-closeup} \end{figure}
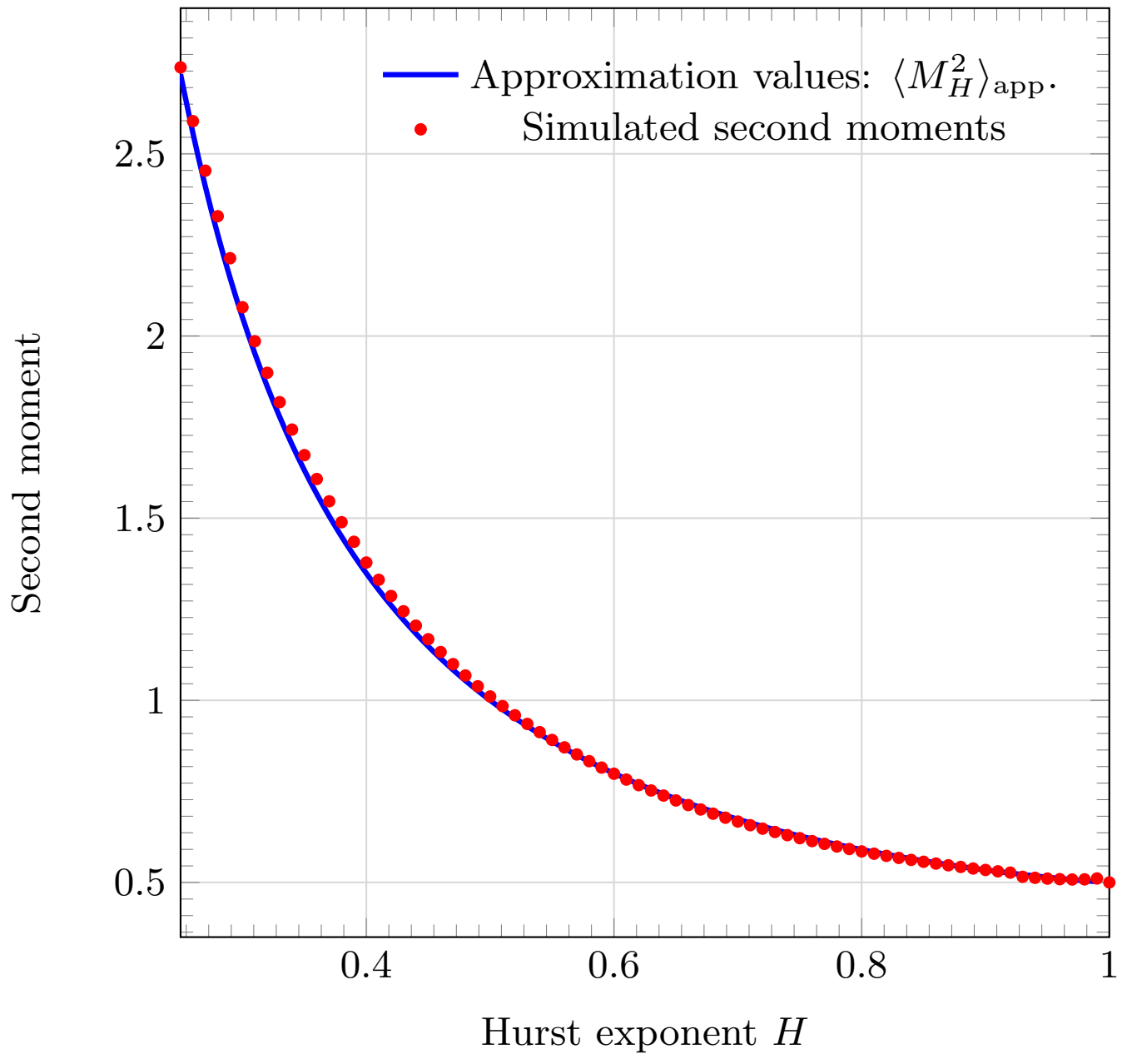

   Figure~\ref{fig:simulation-theory-closeup} compares the simulated second moments of the maxima with the approximation given by Eq.~\eqref{eq:secon_order_appr} over the range \[ 0.25\leq H\leq0.99. \]
These values are strikingly close especially for $H \geq 1/2$. 

 For example, at $H=1/2$ the simulation gave
\begin{equation}
\langle M_{1/2, n}^2\rangle\simeq 1.01056
\end{equation}
close to the exact value $\langle M_{1/2}^2 \rangle = 1$.  At $H=0.6$, the simulation gave
\begin{equation}
\langle M_{0.6,n}^2 \rangle \simeq 0.79849,
\end{equation}
whereas Eq.~\eqref{eq:conjecture} gives $0.79950$.
Also at $H = 0.30$, the simulated value is $2.07901$ while Eq.~\eqref{eq:conjecture} gives $2.05082$. All these values are strikingly close.

\section{Limit of the variance as $H\downarrow0$} 
\label{subsec:discrete-H-zero} Numerical simulations suggest that \[ \operatorname{Var}(M_H)\longrightarrow\frac12 \qquad\text{as }H\downarrow0. \] To gain some insight into this observed behaviour, fix distinct grid points $0<t_1<\cdots<t_n\leq1$. From the fBm covariance, \[ R_H(t_i,t_j) = \frac12\left( t_i^{2H}+t_j^{2H}-|t_i-t_j|^{2H} \right), \] we obtain, as $H\downarrow0$, \[ R_H(t_i,t_j)\longrightarrow \begin{cases} 1, & i=j,\\[1mm] 1/2, & i\neq j. \end{cases} \] 
The limiting Gaussian field therefore admits the representation \begin{equation} X_i=\frac{Z_0+Z_i}{\sqrt2}, \qquad i=1,\ldots,n, \label{eq:H-zero-Gaussian-representation} \end{equation} where $Z_0,Z_1,\ldots,Z_n$ are independent standard normal variables. Let \[ M_n^{(0)}=\max(0,X_1,\ldots,X_n), \qquad Y_n=\max_{1\leq i\leq n}Z_i. \] Then \begin{equation} M_n^{(0)} = \left( \frac{Z_0+Y_n}{\sqrt2} \right)^+. \label{eq:H-zero-maximum-representation} \end{equation} Since $Y_n\to\infty$ almost surely, the positive-part operation becomes asymptotically irrelevant in $L^2$. Hence \[ \operatorname{Var}(M_n^{(0)}) = \frac12+\frac12\operatorname{Var}(Y_n)+o(1). \] Let \[ b_n=\Phi^{-1}\left(1-\frac1n\right) \sim\sqrt{2\log n}. \] The classical Gaussian extreme-value theorem gives \[ b_n(Y_n-b_n) \overset{\mathrm d}{\longrightarrow} \mathcal{D}_G, \] where $\mathcal{D}_G$ has the standard Gumbel distribution. Consequently, \[ \operatorname{Var}(Y_n) \sim \frac{\pi^2}{6b_n^2} \sim \frac{\pi^2}{12\log n} \longrightarrow0 \] \cite{LeadbetterLindgrenRootzen1983}. Therefore, \begin{equation}  \lim_{n\to\infty} \operatorname{Var}(M_n^{(0)}) = \frac12.  \label{eq:H-zero-variance-limit} \end{equation}
Moreover,\[ M_n^{(0)}-\mathbb \langle M_n^{(0)}\rangle = \frac{Z_0}{\sqrt2} + \frac{Y_n-\mathbb \langle Y_n\rangle }{\sqrt2} + R_n, \] where the remainder $R_n$ converges to zero in mean square as $n\to\infty$. 
Since \(\operatorname{Var}(Y_n)\to0\), both the second term and the remainder converge to zero in mean square. Therefore, \begin{equation} M_n^{(0)}-\left\langle M_n^{(0)}\right\rangle \overset{\mathrm d}{\longrightarrow} N\left(0,\frac12\right). \label{eq:H-zero-centered-normal-limit} \end{equation} 
Thus, the independent maximum determines the diverging location, whereas the limiting centered fluctuations arise entirely from the common Gaussian component $Z_0/\sqrt2$. The mean itself diverges but its centered fluctuations remain bounded and converge to a Gaussian distribution with variance $1/2$. 
The probabilistic mechanism is transparent from the decomposition \[ M_n^{(0)} \approx \frac{Y_n}{\sqrt2} + \frac{Z_0}{\sqrt2}. \] The first term determines the diverging location of the maximum. Its centered fluctuations vanish as $n\to\infty$. The second term is a common Gaussian shift, shared by all grid points, and survives in the limit. Its variance is exactly $1/2$, which explains Eq.~\eqref{eq:H-zero-variance-limit}. 
This discrete analysis provides a natural explanation for the numerical observation that, as $H\downarrow0$, $\operatorname{Var}(M_H)$ tends to $1/2$ and the centered maximum tends to a Gaussian distribution with variance $1/2$.

\begin{figure}[t] \centering \includegraphics[ width=0.45\textwidth ]{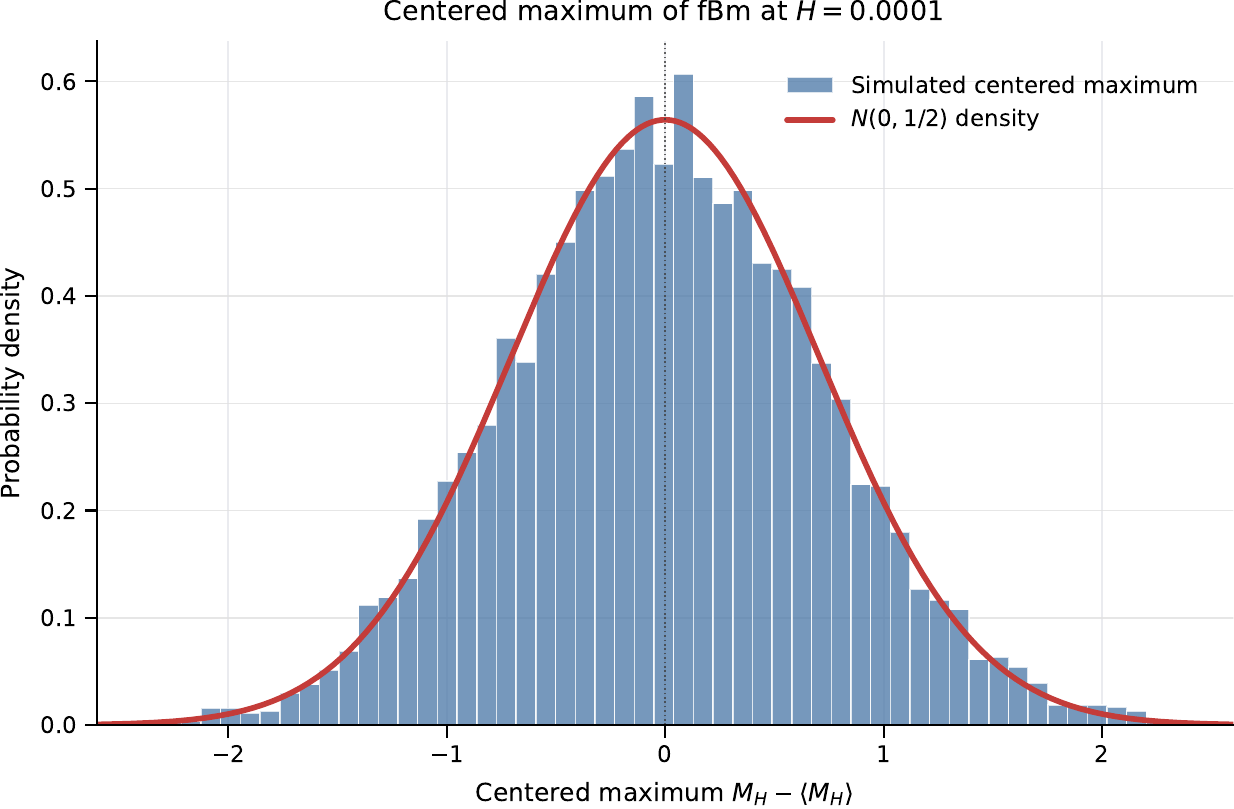} \caption{Density-normalized histogram of the centered maximum $M_H-\langle M_H\rangle$ for $H=0.0001$, obtained from $10\,000$ simulated paths. The solid curve is the density of $N(0,1/2)$.} \label{fig:small-H-histogram} \end{figure}

For $H=10^{-4}$, simulations of $10^4$ paths, using a grid with $2^{20}$ intervals, gave a centered variance of $0.5173$, with skewness $0.0484$ and excess kurtosis $-0.0219$. The histogram in Fig.~\ref{fig:small-H-histogram} agrees closely with the density of $N(0,1/2)$. A bootstrap-calibrated Kolmogorov--Smirnov test gave $p=0.211$, so the Gaussian reference law is not rejected. 

The three limiting regimes suggest an intriguing qualitative picture. After centering, the maximum appears to approach a Gaussian distribution as \(H\downarrow0\), whereas the uncentered maximum is half-normal at \(H=1/2\) and reduces to \(Z^+\), with \(Z\sim N(0,1)\), at \(H=1\). The distributions at intermediate values of \(H\) are not known explicitly. Understanding how their shapes vary with \(H\) remains an open problem.


\section{Conclusion} We have proposed a one-dimensional integral approximation for the second moment of the maximum of fractional Brownian motion. The construction is motivated by a second-moment identity for maxima of partial sums with exchangeable increments, although fractional Gaussian noise is generally not exchangeable. The approximation is exact at \(H=1/2\) and \(H=1\) and agrees closely with the numerical estimates throughout the range investigated. 
Its first-order expansion about \(H=1/2\), however, yields a coefficient that differs slightly from the corresponding correction formally extracted from perturbative results. 
The proposed expression should therefore be regarded as a highly accurate approximation rather than an exact identity. 

Beyond this particular approximation, the results raise a broader structural question concerning the relationship between exchangeability and fractional Brownian motion. 
Although the increments of fractional Brownian motion are not exchangeable, calculations motivated by exchangeability produce remarkably accurate approximations to moments of its maximum, whose exact values remain inaccessible by currently known methods. 
This suggests that the exchangeability-based construction may be capturing a weaker, approximate, or asymptotic symmetry that has not yet been identified. 

The observed connection also suggests a practical approximation strategy. When an extremal quantity of fractional Brownian motion is analytically intractable, one may formally treat the increments as exchangeable, derive the resulting candidate formula, and assess its accuracy through numerical simulation and any available exact or perturbative results. Such a procedure cannot establish an exact identity, but it may yield useful approximations when no exact solution is available.

What feature of fractional Brownian motion is responsible for the success of these formally exchangeable arguments? Can this phenomenon be given a rigorous probabilistic interpretation, and might it lead to systematic approximations for other extremal observables of correlated Gaussian processes? We hope that the results presented here will encourage probabilists and physicists to investigate these questions more closely.

\section*{Data availability}
No external datasets were used. The numerical results were generated from the simulations described in the text.

\appendix

\section{Regularized evaluation of the perturbative integral}
\label{app:integrals}
The integral in Eq.~\eqref{eq:J2-definition} is not straightforward to evaluate directly. We therefore introduce a convergence factor, combine the resulting regularized contributions, and only then remove the regulator by taking the appropriate limit.
 Choose $\lambda >1$ and let    \[J(\lambda) = \int_0^{\infty} \sqrt{\frac{2}{\pi}} e^{-\lambda y^2/2} y^2\,  G(y) dy\] for $\lambda >1$. This  integral is obtained by replacing the quantity $e^{-y^2/2}$ in  Eq.~\eqref{eq:J2-definition}  with $e^{-\lambda y^2/2}$. 

 Write \[J(\lambda) = \sqrt{2/\pi} I(\lambda)\] with 
\begin{equation}
L(\lambda)
=J_{\mathrm{hyp}}(\lambda)
+J_{\mathrm{pol}}(\lambda)
+J_{\mathrm{erfi}}(\lambda)
+J_{\exp}(\lambda)
+J_{\log}(\lambda),
\label{eq:I-decomposition}
\end{equation}
where the five terms correspond, respectively, to the hypergeometric,
polynomial, imaginary-error-function, exponential, and logarithmic
parts of $G$ (in Eq.~\eqref{eq:DW-G}). 

\subsection{The hypergeometric contribution}

The first contribution is
\begin{equation}
J_{\mathrm{hyp}}(\lambda)
=
\frac16\int_0^\infty
y^6 e^{-\lambda y^2/2}
{}_2F_2\left(1,1;\frac52,3;\frac{y^2}{2}\right)dy.
\label{eq:Ihyp-start}
\end{equation}
With $z=y^2/2$, one has $y^6dy=4\sqrt2\,z^{5/2}dz$, and hence
\begin{equation}
J_{\mathrm{hyp}}(\lambda)
=
\frac{2\sqrt2}{3}
\int_0^\infty
z^{5/2}\, e^{-\lambda z}
{}_2F_2\left(1,1;\frac52,3;z\right)dz.
\end{equation}
We use the Laplace transform of the generalized hypergeometric function given in Erd\'elyi et al.\ \cite[p.~220, \S4.23]{ErdelyiEtAl1954}
\begin{equation}
\int_0^\infty t^{\nu-1} e^{-s t} 
{}_pF_q(\bm a;\bm b; \omega t)\,dt = 
\Gamma(\nu) s^{-\nu}
{}_{p+1}F_q(\bm a, \nu;\bm b;\omega s^{-1}),
\qquad 
\end{equation}
to obtain
\begin{equation}
J_{\mathrm{hyp}}(\lambda)
=
\frac{5\sqrt{2\pi}}{4\lambda^{7/2}}
{}_3F_2\left(
1,1,\frac72;
\frac52,3;
\frac1\lambda
\right).
\label{eq:Ihyp-final}
\end{equation}

\subsection{The polynomial contribution}

The polynomial part gives
\begin{equation}
J_{\mathrm{pol}}(\lambda)
=-3\int_0^\infty y^4\,e^{-\lambda y^2/2}dy.
\end{equation}
The standard Gaussian moment formula
\begin{equation}
\int_0^\infty y^m\,e^{-\lambda y^2/2}dy
=
2^{(m-1)/2}\lambda^{-(m+1)/2}
\Gamma\left(\frac{m+1}{2}\right)
\end{equation}
yields
\begin{equation}
J_{\mathrm{pol}}(\lambda)
=-\frac{9\sqrt{2\pi}}{2\lambda^{5/2}}.
\label{eq:Ipol-final}
\end{equation}

\subsection{The imaginary-error-function contribution}
We have 
\[J_{\mathrm{erfi}}(\lambda) = \pi\int_0^\infty y^2 (1-y^2)\operatorname{erfi}\left(\frac y{\sqrt2}\right) e^{-\lambda y^2/2} dy.\]
Set
\begin{equation}
J_0(\lambda)
:=
\int_0^\infty
e^{-\lambda y^2/2}
\operatorname{erfi}\left(\frac{y}{\sqrt2}\right)dy.
\end{equation}
For $\lambda>1$,
\begin{equation}
J_0(\lambda)
=
\sqrt{\frac{2}{\pi \lambda}}
\mbox{arctanh}\left(\frac1{\sqrt \lambda}\right).
\label{eq:J0}
\end{equation}
Differentiation under the integral sign gives
\begin{equation}
\int_0^\infty y^2\,e^{-\lambda y^2/2}
\operatorname{erfi}\left(\frac{y}{\sqrt2}\right)dy
=-2J_0'(\lambda),
\end{equation}
and
\begin{equation}
\int_0^\infty y^4\, e^{-\lambda y^2/2}
\operatorname{erfi}\left(\frac{y}{\sqrt2}\right)dy
=4J_0''(\lambda).
\end{equation}
Therefore
\begin{align}
J_{\mathrm{erfi}}(\lambda) = 
\pi\left(-2J_0'(\lambda)-4J_0''(\lambda)\right).
\end{align}
This implies 
\begin{equation}
\begin{aligned}
J_{\mathrm{erfi}}(\lambda)
=
\sqrt{2\pi}\Bigg(
&\frac{\lambda-3}{\lambda^{5/2}}
\mbox{arctanh}(\lambda^{-1/2})
+
\frac{\lambda^2-6\lambda+3}{\lambda^2(\lambda-1)^2}
\Bigg).
\end{aligned}
\label{eq:Ierfi-final}
\end{equation}

\subsection{The exponential contribution}

The exponential term is
\begin{equation}
J_{\exp}(\lambda) = 
\sqrt{2\pi}
\int_0^\infty y^3
\,e^{-(\lambda-1)y^2/2}dy \\
 =  \frac{2\sqrt{2\pi}}{(\lambda-1)^2}.
\label{eq:Iexp-final}
\end{equation}

\subsection{The logarithmic contribution}

The final contribution is
\begin{equation}
J_{\log}(\lambda)
=
\int_0^\infty
y^2(y^2-2)\,e^{-\lambda y^2/2}
\left[\log(2y^2)+\gamma_{\mathrm E}\right]dy.
\label{eq:Ilog-start}
\end{equation}
Introduce the Mellin moment
\begin{equation} \begin{aligned} M_{\lambda}(p) &:=\int_0^\infty y^p e^{-\lambda y^2/2}\,dy \\ &=2^{(p-1)/2}\lambda^{-(p+1)/2} \Gamma\left(\frac{p+1}{2}\right). \end{aligned} \label{eq:Gaussian-Mellin-moment} \end{equation}
Then
\begin{equation} \begin{aligned} &\int_0^\infty y^p e^{-\lambda y^2/2} \left[\log(2y^2)+\gamma_{\mathrm E}\right]\,dy \\ &\qquad= \left(\log 2+\gamma_{\mathrm E}\right)M_{\lambda}(p) +2\frac{\partial M_{\lambda}(p)}{\partial p}. \end{aligned} \label{eq:Gaussian-log-Mellin-moment} \end{equation}
Applying this identity with $p=4$ and $p=2$ gives
\begin{equation}
J_{\log}(\lambda)
=
\frac{\sqrt{2\pi}}{2\lambda^{5/2}}
\left(
8-4\lambda+(2\lambda-3)\log \lambda
\right).
\label{eq:Ilog-final}
\end{equation}

\subsection{Closed form for the regulated integral}

Combining \eqref{eq:Ihyp-final}, \eqref{eq:Ipol-final},
\eqref{eq:Ierfi-final}, \eqref{eq:Iexp-final}, and
\eqref{eq:Ilog-final}, and using
\begin{equation}
\frac{\lambda^2-6\lambda+3}{\lambda^2(\lambda-1)^2}
+
\frac{2}{(\lambda-1)^2}
=
\frac3{\lambda^2},
\end{equation}
we obtain the following expression.
For $\lambda>1$,
\begin{equation}
\begin{aligned}
L(\lambda)
=
\sqrt{2\pi}\Bigg[
&\frac{5}{4\lambda^{7/2}}
{}_3F_2\left(
1,1,\frac72;
\frac52,3;
\frac1\lambda
\right)
\\
&+
\frac{\lambda-3}{\lambda^{5/2}}
\mbox{arctanh}\left(\frac1{\sqrt \lambda}\right)
+
\frac3{\lambda^2}
\\
&+
\frac{(2\lambda-3)\log \lambda-4\lambda-1}{2\lambda^{5/2}}
\Bigg].
\end{aligned}
\label{eq:Ilambda-closed}
\end{equation}
 
The closed form \eqref{eq:Ilambda-closed} is finite for every $\lambda>1$.
However, its hypergeometric and inverse-hyperbolic terms are separately
singular as $\lambda\downarrow1$. The singular terms must therefore be
expanded and combined before the regulator is removed.

\subsection{Removal of the regulator}
Using the defining hypergeometric series, one obtains the elementary reduction\[ {}_3F_2\left( 1,1,\frac72;\frac52,3;z \right) = \frac{2}{5z^2} \left[z-(3z-1)\log(1-z)\right], \qquad |z|<1. \] Indeed, both sides have the coefficient \[ \frac{2(2n+5)}{5(n+1)(n+2)} \] at order $z^n$.
Setting $z=1/\lambda$ yields 
\begin{equation} {}_3F_2\left( 1,1,\frac72;\frac52,3;\frac1\lambda \right) = \frac25\left[1-2\log(\lambda-1)\right]+o(1) 
\label{eqn:C1}
\end{equation}
as $\lambda\downarrow1$.
Also
\begin{equation}
\mbox{arctanh}\left(\frac1{\sqrt \lambda}\right) = \log 2 - (1/2)\log(\lambda-1)  + o(1). 
\label{eqn:C2}
\end{equation}
Then \eqref{eqn:C1} and \eqref{eqn:C2} in \eqref{eq:Ilambda-closed} implies that 
\[\lim_{\lambda\downarrow1} I(\lambda) = (1-2 \log(2)) \sqrt{2\pi}.\]
Then
$$J = \lim_{\lambda\downarrow1} J(\lambda) = 2(1-2\log(2))$$
Hence 
\begin{equation}
 \left.\frac d{dH}\langle M_H^2 \rangle \right|_{H=1/2}
 =J-2 = -4 \log2.
 \end{equation}
 \newpage

\bibliographystyle{apsrev4-2}

\end{document}